\documentclass[reqno,12pt]{amsart}
\usepackage{amssymb}
\usepackage{amsfonts}
\usepackage{amsmath}
\usepackage{graphicx}
\usepackage{amscd,amsthm}

\newtheorem{theorem}{Theorem}
\theoremstyle{plain}

\newtheorem{corollary}{Corollary}

\newtheorem{lemma}{Lemma}

\newtheorem{remark}{Remark}

\numberwithin{equation}{section}

\input{tcilatex}

\begin{document}
\author{}
\title{}
\maketitle

\begin{center}
\thispagestyle{empty} \pagestyle{myheadings} 
\markboth{\bf Rahime Dere and Yilmaz
Simsek}{\bf Unification of the three families of generalized Apostol type polynomials on the Umbral algebra}

\textbf{\Large Unification of the three families of generalized Apostol type
polynomials on the Umbral algebra}

\bigskip

\textbf{Rahime Dere} \textbf{and Yilmaz Simsek}\\[0pt]

Department of Mathematics, Faculty of Art and Science University of Akdeniz
TR-07058 Antalya, Turkey \\[0pt]

E-mail\textbf{: rahimedere@gmail.com and ysimsek@akdeniz.edu.tr}\\[0pt]

\medskip

\textbf{{\large {Abstract}}}\medskip
\end{center}

\begin{quotation}
The aim of this paper is to investigate and introduce some new identities
related to the unification and generalization of the three families of
generalized Apostol type polynomials,\ which are Apostol-Bernoulli,
Apostol-Euler, and Apostol-Genocchi polynomials, on the modern theory of the
Umbral calculus and algebra. We also introduce some operators. Recently,
Ozden constructed generating function of the unification of the Apostol type
polynomials (see Ozden [H. Ozden, \textit{AIP Conf. Proc.} 1281, (2010),
1125-1227.]). By using this generating function, we derive many properties
of these polynomials. We give relations between these polynomials and
Stirling numbers.
\end{quotation}

\bigskip


\bigskip


\section{Introduction, Definitions and Preliminaries}

Throughout of this paper, we use the following notations:%
\begin{equation*}
\mathbb{N}:=\left\{ 1,2,3,\cdots \right\} \text{ and }\mathbb{N}_{0}:=%
\mathbb{N\cup }\left\{ 0\right\} \text{.}
\end{equation*}%
Here $\mathbb{Z}$ denotes the set of integers, $\mathbb{R}$ denotes the set
of real numbers, $\mathbb{R}^{+\text{ }}$denotes the set of positive real
numbers and $\mathbb{C}$ denotes the set of complex numbers. We also assume
that $\ln \left( z\right) $ denotes the \textit{principal branch of the
many-valued function }$\ln \left( z\right) $ with the imaginary part $\Im
\left( \ln \left( z\right) \right) $ constrained by $-\pi <\Im \left( \ln
\left( z\right) \right) \leqq \pi $.%
\begin{equation*}
\binom{m}{n}=\frac{m!}{n!(m-n)!},
\end{equation*}%
\begin{equation*}
\binom{m}{a,b,c}=\frac{m!}{a!b!c!}\text{ and }a+b+c=m,
\end{equation*}%
\begin{equation*}
0^{j}=\left\{ 
\begin{array}{c}
0\text{ if }j\neq 0 \\ 
1\text{ if }j=0,%
\end{array}%
\right.
\end{equation*}%
and

\begin{equation*}
\left( x\right) _{b}=x\left( x-1\right) ...\left( x-b+1\right) ,
\end{equation*}%
where $b\in \mathbb{N}$ and $\left( x\right) _{0}=1$.

Many generating functions of the Bernoulli, Euler and Genocchi type
polynomials have been introduced by many authors. Recently, Ozden \cite%
{Ozden}-\cite{OzdenICNAAM2011} constructed generating functions of the
unification of the Apostol type polynomials $\mathcal{Y}_{n,\beta
}^{(v)}(x;k,a,b)$. These generating functions were investigated rather
systematically by Ozden \textit{et al}. \cite{ozdensimseksrivstava}.

The following generating function of the\textit{\ unification }of the
Apostol type Bernoulli, Euler and Genocchi polynomials $\mathcal{Y}_{n,\beta
}^{(v)}(x;k,a,b)$\ of higher order, which was recently defined by Ozden \cite%
{OzdenICNAAM2011}:%
\begin{equation}
\left( \frac{2^{1-k}t^{k}}{\beta ^{b}e^{t}-a^{b}}\right)
^{v}e^{xt}=\dsum\limits_{n=0}^{\infty }\mathcal{Y}_{n,\beta }^{(v)}(x;k,a,b)%
\frac{t^{n}}{n!},  \label{b1}
\end{equation}%
where%
\begin{equation*}
\left\vert t+b\ln \left( \frac{\beta }{a}\right) \right\vert <2\pi ;\text{ }%
x\in 
\mathbb{R}
\end{equation*}%
and $k\in 
\mathbb{N}
_{0}$; $v\in 
\mathbb{N}
$; $a,b\in 
\mathbb{R}
^{+}$; $\beta \in \mathbb{C}$.

Observe that%
\begin{equation}
\mathcal{Y}_{n,\beta }^{(1)}(x;k,a,b)=\mathcal{Y}_{n,\beta }(x;k,a,b)\text{ (%
}n\in 
\mathbb{N}
\text{)}  \label{ba1}
\end{equation}%
which is defined by Ozden \cite{Ozden}. Recently, Ozden \textit{et al}. \cite%
{ozdensimseksrivstava} introduced many properties of these polynomials.

\begin{remark}
By substituting $x=0$ in the generating function (\ref{b1}), we obtain the
corresponding unification and generalization of the generating functions of
Bernoulli, Euler and Genocchi numbers of higher order. Thus, we have%
\begin{equation*}
\mathcal{Y}_{n,\beta }^{(v)}(0;k,a,b)=\mathcal{Y}_{n,\beta }^{(v)}(k,a,b).
\end{equation*}
\end{remark}

\begin{remark}
By setting $b=v=1$ and $\beta =1$ into (\ref{ba1}), we have $\mathcal{Y}%
_{n,1}(x;k,a,1)=D_{n}(x;a,k)$ (see, for details, \cite{Karande}).
\end{remark}

\begin{remark}
By substituting $a=b=k=1$ into (\ref{b1}), one has the Apostol-Bernoulli
polynomials $\mathcal{Y}_{n,\beta }^{(v)}(x;1,1,1)=\mathcal{B}%
_{n}^{(v)}(x,\beta ),$which are defined by means of the following generating
function:%
\begin{equation*}
\left( \frac{t}{\beta e^{t}-1}\right) ^{v}e^{xt}=\sum_{n=0}^{\infty }%
\mathcal{B}_{n}^{(v)}(x,\beta )\frac{t^{n}}{n!},\text{ }(\left\vert t+\log
\beta \right\vert <2\pi )
\end{equation*}%
(see, for details, \cite{ozdenSimsekAML}, \cite{Ozden}, \cite%
{ozdensimseksrivstava}, \cite{LuoSrivastava}, \cite{simsekKimDkim}, \cite%
{simsekASCM2008}, \cite{SIMSEKjnt-2005}, \cite{simtwistedhqbernoullijmaa}, 
\cite{sirKIMsim}; see also the references cited in each of these earlier
works).
\end{remark}

\begin{remark}
If we substitute $b=v=1$, $k=0$ and$\ a=-1$ into (\ref{b1}), \ then we have
the Apostol-Euler polynomials $\mathcal{Y}_{n,\beta }(x;0,-1,1)=\mathcal{E}%
_{n}(x,\beta )$, which is defined by the following generating function:%
\begin{equation*}
\frac{2e^{xt}}{\beta e^{t}+1}=\sum_{n=0}^{\infty }\mathcal{E}_{n}(x,\beta )%
\frac{t^{n}}{n!},\text{ }(\left\vert t+\log \beta \right\vert <\pi )
\end{equation*}%
(see, for details, \cite{ozdenSimsekAML}, \cite{Ozden}, \cite%
{ozdensimseksrivstava}, \cite{LuoSrivastava}, \cite{simsekKimDkim}, \cite%
{simsekASCM2008}, \cite{SIMSEKjnt-2005}, \cite{sirKIMsim}; see also the
references cited in each of these earlier works).
\end{remark}

\begin{remark}
By substituting $b=1=v$, $k=1$ and$\ a=-1$ into (\ref{b1}), one has the
Apostol-Genocchi polynomials $\mathcal{Y}_{n,\beta }(x;1,-1,1)=\frac{1}{2}%
\mathcal{G}_{n}(x,\beta )$, which is defined by the following generating
function:%
\begin{equation*}
\frac{2te^{xt}}{\beta e^{t}+1}=\sum_{n=0}^{\infty }\mathcal{G}_{n}(x,\beta )%
\frac{t^{n}}{n!},\text{ }(\left\vert t+\log \beta \right\vert <\pi )
\end{equation*}%
(see, for details, \cite{Luo3}, \cite{ozdenSimsekAML}, \cite{Ozden}, \cite%
{ozdensimseksrivstava}, \cite{Jang}; see also the references cited in each
of these earlier works).
\end{remark}

\begin{remark}
Substituting $\beta =b=$ $k=a=v=1$ into (\ref{b1}), we have $\mathcal{Y}%
_{n,1}(x;1,1,1)=B_{n}(x)$, where $B_{n}(x)$\ denotes the classical Bernoulli
polynomials, which are defined by the following generating function:%
\begin{equation*}
\frac{te^{xt}}{e^{t}-1}=\sum_{n=0}^{\infty }B_{n}(x)\frac{t^{n}}{n!},\text{ }%
(\left\vert t\right\vert <2\pi )
\end{equation*}%
(see, for example, \cite{blasiack}-\cite{sirKIMsim}; see also the references
cited in each of these earlier works).
\end{remark}

\begin{remark}
Likely, substituting $\beta =b=$ $k=v=1$, and $a=-1$ into (\ref{b1}), we have%
\begin{equation*}
\mathcal{Y}_{n,1}(x;0,-1,1)=E_{n}(x),
\end{equation*}%
where $E_{n}(x)$\ denotes the classical Euler polynomials, which are defined
by the following generating function:%
\begin{equation*}
\frac{2e^{xt}}{e^{t}+1}=\sum_{n=0}^{\infty }E_{n}(x)\frac{t^{n}}{n!},\text{ }%
(\left\vert t\right\vert <\pi )
\end{equation*}%
(see, for example, \cite{blasiack}-\cite{sirKIMsim}; see also the references
cited in each of these earlier works).
\end{remark}

We use the following relations and identities of the umbral algebra and
calculus which are given by Roman \cite{Roman}:

Let $P$ be the algebra of polynomials in the single variable $x$ over the
field complex numbers. Let $P^{\ast }$ be the vector space of all linear
functionals on $P$. Let $\left\langle L\mid p(x)\right\rangle $ be the
action of a linear functional $L$ on a polynomial $p(x)$. Let $\tciFourier $
denotes the algebra of formal power series%
\begin{equation*}
f\left( t\right) =\dsum\limits_{k=0}^{\infty }\frac{a_{_{k}}}{k!}t^{k}.
\end{equation*}

This kind of algebra is called an umbral algebra. Each $f\in \tciFourier $
defines a linear functional on $P$ and for all $k\geqslant 0$, $%
a_{k}=\left\langle f\left( t\right) \mid x^{k}\right\rangle $. The order $%
o\left( f\left( t\right) \right) $ of a power series $f\left( t\right) $ is
the smallest integer $k$ for which the coefficient of $t^{k}$ does not
vanish. A series $f\left( t\right) $\ for which $o\left( f\left( t\right)
\right) =1$ is called a delta series. And a series $f\left( t\right) $\ for
which $o\left( f\left( t\right) \right) =0$ is called a invertible series.

Let $f(t)$, $g(t)$ be in $\tciFourier $, we have%
\begin{equation}
\left\langle f(t)g(t)\mid p\left( x\right) \right\rangle =\left\langle
f(t)\mid g(t)p\left( x\right) \right\rangle .  \label{a4}
\end{equation}%
For all $p\left( x\right) $ in $P$, we have%
\begin{equation}
\left\langle e^{yt}\mid p\left( x\right) \right\rangle =p\left( y\right) .
\label{a5}
\end{equation}

The Sheffer polynomials are defined by means of the following generating
function%
\begin{equation}
\dsum\limits_{k=0}^{\infty }\frac{s_{k}\left( x\right) }{k!}t^{k}=\frac{1}{%
g(t)}e^{xt}.  \label{a8}
\end{equation}%
(see, for details, \cite{Roman}; and see also \cite{blasiack}, \cite{Dattoli}%
).

\begin{theorem}
(See, for details, Roman \cite[p. 20, Theorem 2.3.6]{Roman}) Let $f\left(
t\right) $ be a delta series and let $g\left( t\right) $ be an invertible
series. Then there exist a unique sequence $s_{n}\left( x\right) $ of
polynomials satisfying the orthogonality conditions%
\begin{equation}
\left\langle g(t)f(t)^{k}\mid s_{n}(x)\right\rangle =n!\delta _{n,k}
\label{a9}
\end{equation}%
for all $n,k\geq 0$.
\end{theorem}

\begin{remark}
The sequence $s_{n}(x)$\ in (\ref{a9}) is the Sheffer polynomials for pair $%
(g(t),f(t))$, where $g(t)$ must be invertible and $f(t)$\ must be delta
series. The Sheffer polynomials for pair $(g(t),t)$ is the Appell
polynomials or the Appell sequences for $g(t)$. The Bernoulli polynomials,
the Euler polynomials, the Genocchi polynomials, and the polynomials $%
\mathcal{Y}_{n,\beta }^{(v)}(x;k,a,b)$\ are the Appell polynomials cf. \cite%
{blasiack}-\cite{Tempesta}.
\end{remark}

Let%
\begin{equation}
s_{n}\left( x\right) =g(t)^{-1}x^{n},  \label{a10}
\end{equation}%
derivative formula%
\begin{equation}
ts_{n}\left( x\right) =s_{n}^{^{\prime }}\left( x\right) =ns_{n-1}\left(
x\right) ,  \label{a11}
\end{equation}%
recurrence formula%
\begin{equation}
s_{n+1}\left( x\right) =\left( x-\frac{g^{^{\prime }}\left( t\right) }{%
g\left( t\right) }\right) s_{n}\left( x\right) ,  \label{a12}
\end{equation}%
and multiplication formula, for $\alpha \neq 0$%
\begin{equation}
s_{n}\left( \alpha x\right) =\alpha ^{n}\frac{g\left( t\right) }{g\left( 
\frac{t}{\alpha }\right) }s_{n}\left( x\right) .  \label{a14}
\end{equation}%
Proof of equations (\ref{a10})-(\ref{a14}) are given by Roman \cite{Roman}.

\section{The operator $\frac{1}{t}$}

In this section we introduce $\frac{1}{f\left( t\right) }$ operator on the
Umbral algebra. By applying this operator to the Sheffer sequences, we
obtain new identities related to the family of the Sheffer sequences.

\begin{theorem}
\label{Teo-0} Let $n\in 
\mathbb{N}
_{0}$. Let $s_{n}(x)$ be Sheffer sequence for $\left( g\left( t\right)
,f\left( t\right) \right) $. The following relationship holds true:%
\begin{equation*}
\frac{1}{f\left( t\right) }s_{n}\left( x\right) =\frac{1}{n+1}s_{n+1}\left(
x\right) .
\end{equation*}
\end{theorem}

\begin{proof}
We set%
\begin{equation*}
\left\langle g(t)f(t)^{k}\mid \frac{1}{f\left( t\right) }s_{n}(x)\right%
\rangle =\left\langle g(t)f(t)^{k-1}\mid s_{n}(x)\right\rangle .
\end{equation*}%
By using (\ref{a9}) in the above relation, we obtain%
\begin{equation*}
\left\langle g(t)f(t)^{k}\mid \frac{1}{f\left( t\right) }s_{n}(x)\right%
\rangle =n!\delta _{n,k-1}.
\end{equation*}%
From the definition of the Kronecker $\delta $, we get%
\begin{equation*}
\left\langle g(t)f(t)^{k}\mid \frac{1}{f\left( t\right) }s_{n}(x)\right%
\rangle =\left\langle g(t)f(t)^{k}\mid \frac{1}{n+1}s_{n+1}(x)\right\rangle .
\end{equation*}%
After some calculation in the above equation, we obtain the desired result.
\end{proof}

Taking $f\left( t\right) =t$ in Theorem \ref{Teo-0}, we deduce the following
corollary:

\begin{corollary}
Let $n\in 
\mathbb{N}
_{0}$. Let $s_{n}(x)$ be the Appell sequence for $g\left( t\right) $. The
following relationship holds true:%
\begin{equation}
\frac{1}{t}s_{n}\left( x\right) =\frac{1}{n+1}s_{n+1}\left( x\right) .
\label{brd}
\end{equation}
\end{corollary}

\section{Some new identities and relations of the unification of the
Bernoulli, Euler and Genocchi polynomials $\mathcal{Y}_{n,\protect\beta %
}^{(v)}(x;k,a,b)$ of higher order}

In this section, by using properties of the Sheffer sequences and the Appell
sequences, we give some important properties of the polynomials $\mathcal{Y}%
_{n,\beta }^{(v)}(x,k,a,b)$ and the proof of this properties.

By using (\ref{a10})\ and (\ref{b1}), we arrive at the following lemma:

\begin{lemma}
\label{Lemma1} Let $n\in 
\mathbb{N}
_{0}$. The following relationship holds true:%
\begin{equation*}
\mathcal{Y}_{n,\beta }^{(v)}(x;k,a,b)=\left( \frac{2^{1-k}t^{k}}{\beta
^{b}e^{t}-a^{b}}\right) ^{v}x^{n}.
\end{equation*}
\end{lemma}

\begin{lemma}
\label{Lemma1A} Let $n\in 
\mathbb{N}
_{0}$. The following relationship holds true:%
\begin{equation*}
\left\langle (\beta ^{b}e^{t}-a^{b})^{j}\mid \mathcal{Y}_{n,\beta
}(x;k,a,b)\right\rangle =\dsum\limits_{m=0}^{j}\binom{j}{m}\left( -a\right)
^{b\left( j-m\right) }\beta ^{bm}\mathcal{Y}_{n,\beta }(m;k,a,b).
\end{equation*}
\end{lemma}

\begin{proof}
\begin{eqnarray*}
&&\left\langle (\beta ^{b}e^{t}-a^{b})^{j}\mid \mathcal{Y}_{n,\beta
}(x;k,a,b)\right\rangle  \\
&=&\left\langle \dsum\limits_{m=0}^{j}\binom{j}{m}\beta ^{bm}e^{mt}\left(
-a\right) ^{b\left( j-m\right) }\mid \mathcal{Y}_{n,\beta
}(x;k,a,b)\right\rangle  \\
&=&\dsum\limits_{m=0}^{j}\binom{j}{m}\beta ^{bm}\left( -a\right) ^{b\left(
j-m\right) }\left\langle e^{mt}\mid \mathcal{Y}_{n,\beta
}(x;k,a,b)\right\rangle .
\end{eqnarray*}%
Using (\ref{a5}) in this equation, we complete proof of this Lemma.
\end{proof}

\begin{remark}
By substituting $\beta =a=b=k=1$ into Lemma \ref{Lemma1A}, we arrive at a
special case:%
\begin{eqnarray*}
&&\left\langle (e^{t}-1)^{j}\mid \mathcal{Y}_{n,1}(x;1,1,1)\right\rangle
=\left\langle (e^{t}-1)^{j}\mid B_{n}(x)\right\rangle \\
&=&\dsum\limits_{m=0}^{j}\binom{j}{m}\left( -1\right) ^{j-m}B_{n}\left(
m\right) .
\end{eqnarray*}
\end{remark}

\begin{remark}
If we set $\beta =b=1$, $k=0$ and $a=-1$ in Lemma \ref{Lemma1A}, we arrive
at the following result:%
\begin{eqnarray*}
&&\left\langle (e^{t}+1)^{j}\mid \mathcal{Y}_{n,1}(x;0,-1,1)\right\rangle
=\left\langle (e^{t}+1)^{j}\mid E_{n}(x)\right\rangle \\
&=&\dsum\limits_{m=0}^{j}\binom{j}{m}E_{n}\left( m\right) .
\end{eqnarray*}
\end{remark}

\begin{remark}
Substituting $\beta =b=k=1$, and $a=-1$ into Lemma \ref{Lemma1A}, we obtain%
\begin{eqnarray*}
&&\left\langle (e^{t}+1)^{j}\mid \mathcal{Y}_{n,1}(x;1,-1,1)\right\rangle
=\left\langle (e^{t}+1)^{j}\mid \frac{1}{2}G_{n}(x)\right\rangle \\
&=&\frac{1}{2}\dsum\limits_{m=0}^{j}\binom{j}{m}G_{n}\left( m\right) .
\end{eqnarray*}
\end{remark}

\begin{remark}
For $k=a=b=1$ into Lemma \ref{Lemma1A}, we arrive at the following result.%
\begin{eqnarray*}
&&\left\langle (\beta e^{t}-1)^{j}\mid \mathcal{Y}_{n,\beta
}(x;1,1,1)\right\rangle =\left\langle (\beta e^{t}-1)^{j}\mid \mathcal{B}%
_{n}\left( x\right) \right\rangle \\
&=&\dsum\limits_{m=0}^{j}\binom{j}{m}\left( -1\right) ^{j-m}\beta ^{m}%
\mathcal{B}_{n}\left( m,\beta \right) .
\end{eqnarray*}
\end{remark}

\begin{remark}
By putting $k=0$, $b=1$, and $a=-1$ in Lemma \ref{Lemma1A}, we have%
\begin{eqnarray*}
&&\left\langle (\beta e^{t}+1)^{j}\mid \mathcal{Y}_{n,\beta
}(x;0,-1,1)\right\rangle =\left\langle (\beta e^{t}+1)^{j}\mid \mathcal{E}%
_{n}\left( x\right) \right\rangle \\
&=&\dsum\limits_{m=0}^{j}\binom{j}{m}\beta ^{m}\mathcal{E}_{n}\left( m,\beta
\right) .
\end{eqnarray*}
\end{remark}

\begin{remark}
In a special case when $k=$ $b=1$, and $a=-1$, Lemma \ref{Lemma1A} yields:%
\begin{eqnarray*}
&&\left\langle (\beta e^{t}+1)^{j}\mid \mathcal{Y}_{n,\beta
}(x;1,-1,1)\right\rangle =\left\langle (\beta e^{t}+1)^{j}\mid \frac{1}{2}%
\mathcal{G}_{n}\left( x\right) \right\rangle  \\
&=&\frac{1}{2}\dsum\limits_{m=0}^{j}\binom{j}{m}\beta ^{m}\mathcal{G}%
_{n}\left( m,\beta \right) .
\end{eqnarray*}
\end{remark}

\begin{lemma}
\label{Lemma1.a} Then the following identity holds:%
\begin{equation*}
S\left( n,l\right) =\frac{1}{l!}\left\langle \left( e^{t}-1\right) ^{l}\mid
x^{n}\right\rangle
\end{equation*}%
where $S\left( n,l\right) $ is the Stirling numbers of the second kind.
\end{lemma}

\begin{theorem}
\label{Teo-1} If $a\neq \beta $, then we have%
\begin{eqnarray*}
&&\left\langle (\beta ^{b}e^{t}-a^{b})^{j}\mid \mathcal{Y}_{n,\beta
}(x;k,a,b)\right\rangle \\
&=&\frac{\beta ^{b\left( j-1\right) }k!}{2^{k-1}}\binom{n}{k}%
\dsum\limits_{l=0}^{j-1}\frac{\left( j-1\right) !}{\left( j-l-1\right) !}%
\left( 1-\frac{a^{b}}{\beta ^{b}}\right) ^{j-l-1}S\left( n-k,l\right) .
\end{eqnarray*}%
If $a=\beta $, then we have%
\begin{equation*}
\left\langle (e^{t}-1)^{j}\mid \mathcal{Y}_{n,a}(x;k,a,b)\right\rangle =%
\frac{k!\left( j-1\right) !}{2^{k-1}a^{b}}\binom{n}{k}S\left( n-k,j-1\right)
.
\end{equation*}%
where $S(u,v)$ denote the Stirling numbers of the second kind.
\end{theorem}

\begin{proof}
By applying Lemma \ref{Lemma1} to $\left\langle (\beta
^{b}e^{t}-a^{b})^{j}\mid \mathcal{Y}_{n,\beta }(x,k,a,b)\right\rangle $, we
obtain%
\begin{equation*}
\left\langle (\beta ^{b}e^{t}-a^{b})^{j}\mid \mathcal{Y}_{n,\beta
}(x,k,a,b)\right\rangle =\left\langle (\beta ^{b}e^{t}-a^{b})^{j}\mid \frac{%
2^{1-k}t^{k}}{\beta ^{b}e^{t}-a^{b}}x^{n}\right\rangle .
\end{equation*}%
By using (\ref{a4}) and (\ref{a11}) in the above equation, we get%
\begin{equation*}
\left\langle (\beta ^{b}e^{t}-a^{b})^{j}\mid \mathcal{Y}_{n,\beta
}(x;k,a,b)\right\rangle =2^{1-k}\left\langle (\beta
^{b}e^{t}-a^{b})^{j-1}\mid t^{k}x^{n}\right\rangle .
\end{equation*}%
Applying Lemma \ref{Lemma2} in this equation, we obtain%
\begin{equation*}
\left\langle (\beta ^{b}e^{t}-a^{b})^{j}\mid \mathcal{Y}_{n,\beta
}(x;k,a,b)\right\rangle =2^{1-k}\left\langle (\beta
^{b}e^{t}-a^{b})^{j-1}\mid \left( n\right) _{k}x^{n-k}\right\rangle .
\end{equation*}

Hence we get%
\begin{eqnarray*}
&&\left\langle (\beta ^{b}e^{t}-a^{b})^{j}\mid \mathcal{Y}_{n,\beta
}(x;k,a,b)\right\rangle  \\
&=&2^{1-k}\left( n\right) _{k+1}\beta ^{b\left( j-1\right)
}\dsum\limits_{l=0}^{j-1}\binom{j-1}{l}\left( 1-\frac{a^{b}}{\beta ^{b}}%
\right) ^{j-l-1}\left\langle \left( e^{t}-1\right) ^{l}\mid
x^{n-k}\right\rangle .
\end{eqnarray*}%
Thus by Lemma \ref{Lemma1.a}, we obtain%
\begin{eqnarray*}
&&\left\langle (\beta ^{b}e^{t}-a^{b})^{j}\mid \mathcal{Y}_{n,\beta
}(x;k,a,b)\right\rangle  \\
&=&2^{1-k}\left( n\right) _{k+1}\beta ^{b\left( j-1\right)
}\dsum\limits_{l=0}^{j-1}\frac{\left( j-1\right) !}{\left( j-l-1\right) !}%
\left( 1-\frac{a^{b}}{\beta ^{b}}\right) ^{j-l-1}S\left( n-k,l\right) .
\end{eqnarray*}%
After elementary manipulations in the above equation, we complete proof of\
the theorem.
\end{proof}

By combining Lemma \ref{Lemma1A} and Theorem \ref{Teo-1}, we arrive at the
following corollary:

\begin{corollary}
\label{Corollary1} The following relationship holds true:%
\begin{eqnarray*}
&&\dsum\limits_{m=0}^{j}\binom{j}{m}\left( -a\right) ^{b\left( j-m\right)
}\beta ^{bm}\mathcal{Y}_{n,\beta }(x;k,a,b) \\
&=&\frac{\beta ^{b\left( j-1\right) }k!}{2^{k-1}}\binom{n}{k}%
\dsum\limits_{l=0}^{j-1}\frac{\left( j-1\right) !}{\left( j-l-1\right) !}%
\left( 1-\frac{a^{b}}{\beta ^{b}}\right) ^{j-l-1}S\left( n-k,l\right) .
\end{eqnarray*}
\end{corollary}

\begin{remark}
Theorem \ref{Teo-1} provides us with a generalized and unification of the
linear operator $\left\langle .\mid .\right\rangle $ related to the Apostol
type polynomials.
\end{remark}

\begin{remark}
By setting $\beta =a=b=k=1$ in Theorem \ref{Teo-1}, we arrive at the
well-known relation which was proved by Roman \cite[p. 94]{Roman}:%
\begin{eqnarray*}
&&\left\langle (e^{t}-1)^{j}\mid \mathcal{Y}_{n,1}(x;1,1,1)\right\rangle
=\left\langle (e^{t}-1)^{j}\mid B_{n}(x)\right\rangle \\
&=&n\left( j-1\right) !S\left( n-1,j-1\right) .
\end{eqnarray*}%
And using Corollary \ref{Corollary1}, we get%
\begin{equation*}
\dsum\limits_{m=0}^{j}\binom{j}{m}\left( -1\right) ^{j-m}B_{n}\left(
m\right) =n\left( j-1\right) !S\left( n-1,j-1\right) .
\end{equation*}
\end{remark}

\begin{remark}
By putting $\beta =b=k=1$ and $a=-1$ into Theorem \ref{Teo-1}, we have%
\begin{equation*}
\left\langle (e^{t}+1)^{j}\mid \mathcal{Y}_{n,1}(x;1,-1,1)\right\rangle
=\left\langle (e^{t}+1)^{j}\mid \frac{1}{2}G_{n}(x)\right\rangle .
\end{equation*}%
From this equation, we get%
\begin{equation*}
\left\langle (e^{t}+1)^{j}\mid G_{n}(x)\right\rangle =n\left( j-1\right)
!\dsum\limits_{l=0}^{j-1}\frac{2^{j-l}}{\left( j-l-1\right) !}S\left(
n-1,l\right)
\end{equation*}%
cf. \cite[Theorem 2]{DereSimsek}. And using Corollary \ref{Corollary1}, we
have%
\begin{equation*}
\dsum\limits_{m=0}^{j}\binom{j}{m}G_{n}\left( m\right) =2n\left( j-1\right)
!\dsum\limits_{l=0}^{j-1}\frac{2^{j-l}}{\left( j-l-1\right) !}S\left(
n-1,l\right) .
\end{equation*}
\end{remark}

\begin{remark}
Upon substituting $\beta =b=1$, $k=0$ and $a=-1$ into Theorem \ref{Teo-1},
we obtain%
\begin{eqnarray*}
&&\left\langle (e^{t}+1)^{j}\mid \mathcal{Y}_{n,1}(x;0,-1,1)\right\rangle
=\left\langle (e^{t}+1)^{j}\mid E_{n}(x)\right\rangle \\
&=&\dsum\limits_{l=0}^{j-1}\frac{\left( j-1\right) !}{\left( j-l-1\right) !}%
2^{j-l}S\left( n,l\right) .
\end{eqnarray*}%
And using Corollary \ref{Corollary1}, we have%
\begin{equation*}
\dsum\limits_{m=0}^{j}\binom{j}{m}E_{n}\left( m\right)
=\dsum\limits_{l=0}^{j-1}\frac{\left( j-1\right) !}{\left( j-l-1\right) !}%
2^{j-l}S\left( n,l\right) .
\end{equation*}
\end{remark}

\begin{remark}
If we substitute $a=b=k=1$ into Theorem \ref{Teo-1}, we get a special case
of the generalized Bernoulli polynomials $\mathcal{Y}_{n,\beta }(x;k,a,b)$,
that is, the so-called Apostol-Bernoulli polynomials $\mathcal{B}_{n}\left(
x,\beta \right) $:%
\begin{eqnarray*}
&&\left\langle (\beta e^{t}-1)^{j}\mid \mathcal{Y}_{n,\beta
}(x;1,1,1)\right\rangle =\left\langle (\beta e^{t}-1)^{j}\mid \mathcal{B}%
_{n}\left( x,\beta \right) \right\rangle \\
&=&n\beta ^{j-1}\dsum\limits_{l=0}^{j-1}\frac{\left( j-1\right) !}{\left(
j-l-1\right) !}\left( 1-\frac{1}{\beta }\right) ^{j-l-1}S\left( n-1,l\right)
.
\end{eqnarray*}%
And using Corollary \ref{Corollary1}, we have%
\begin{equation*}
\dsum\limits_{m=0}^{j}\binom{j}{m}\left( -1\right) ^{j-m}\beta ^{m}\mathcal{B%
}_{n}\left( m,\beta \right) =n\beta ^{j-1}\dsum\limits_{l=0}^{j-1}\frac{%
\left( j-1\right) !}{\left( j-l-1\right) !}\left( 1-\frac{1}{\beta }\right)
^{j-l-1}S\left( n-1,l\right) .
\end{equation*}
\end{remark}

\begin{remark}
By putting $b=$ $k=1$, and $a=-1$ in Theorem \ref{Teo-1}, we get a special
case of the generalized Bernoulli polynomials $\mathcal{Y}_{n,\beta
}(x;k,a,b)$, that is, the so-called Apostol-Genocchi polynomials $\mathcal{G}%
_{n}\left( x,\beta \right) :$%
\begin{eqnarray*}
&&\left\langle (\beta e^{t}+1)^{j}\mid \mathcal{Y}_{n,\beta
}(x;1,-1,1)\right\rangle =\left\langle (\beta e^{t}+1)^{j}\mid \frac{1}{2}%
\mathcal{G}_{n}\left( x,\beta \right) \right\rangle \\
&=&n\beta ^{j-1}\dsum\limits_{l=0}^{j-1}\frac{\left( j-1\right) !}{\left(
j-l-1\right) !}\left( 1+\frac{1}{\beta }\right) ^{j-l-1}S\left( n-1,l\right)
.
\end{eqnarray*}%
And using Corollary \ref{Corollary1}, we have%
\begin{equation*}
\dsum\limits_{m=0}^{j}\binom{j}{m}\beta ^{m}\mathcal{G}_{n}\left( m,\beta
\right) =2n\beta ^{j-1}\dsum\limits_{l=0}^{j-1}\frac{\left( j-1\right) !}{%
\left( j-l-1\right) !}\left( 1+\frac{1}{\beta }\right) ^{j-l-1}S\left(
n-1,l\right) .
\end{equation*}
\end{remark}

\begin{remark}
Taking $b=1$, $k=0$, and $a=-1$ into Theorem \ref{Teo-1}, we get a special
case of the generalized Bernoulli polynomials $\mathcal{Y}_{n,\beta
}(x;k,a,b)$, that is, the so-called Apostol-Euler polynomials $\mathcal{E}%
_{n}\left( x,\beta \right) :$%
\begin{eqnarray*}
&&\left\langle (\beta e^{t}+1)^{j}\mid \mathcal{Y}_{n,\beta
}(x;0,-1,1)\right\rangle =\left\langle (\beta e^{t}+1)^{j}\mid \mathcal{E}%
_{n}\left( x,\beta \right) \right\rangle \\
&=&2\beta ^{j-1}\dsum\limits_{l=0}^{j-1}\frac{\left( j-1\right) !}{\left(
j-l-1\right) !}\left( 1+\frac{1}{\beta }\right) ^{j-l-1}S\left( n-1,l\right)
.
\end{eqnarray*}%
And using Corollary \ref{Corollary1}, we have%
\begin{equation*}
\dsum\limits_{m=0}^{j}\binom{j}{m}\beta ^{m}\mathcal{E}_{n}\left( m,\beta
\right) =2\beta ^{\left( j-1\right) }\dsum\limits_{l=0}^{j-1}\frac{\left(
j-1\right) !}{\left( j-l-1\right) !}\left( 1+\frac{1}{\beta }\right)
^{j-l-1}S\left( n,l\right) .
\end{equation*}
\end{remark}

By using (\ref{a11}), we arrive at the following lemma:

\begin{lemma}
\label{Lemma2} We have%
\begin{equation*}
t\mathcal{Y}_{n,\beta }^{(v)}(x;k,a,b)=n\mathcal{Y}_{n-1,\beta
}^{(v)}(x;k,a,b).
\end{equation*}
\end{lemma}

\begin{remark}
\label{Remark1}A second proof of Lemma \ref{Lemma2} is also obtained from (%
\ref{b1}) by using derivative with respect to $x$.
\end{remark}

By applying (\ref{brd}) to the polynomials $\mathcal{Y}_{n,\beta
}^{(v)}(x;k,a,b)$, we arrive at the following Lemma:

\begin{theorem}
The following integral operator holds true:%
\begin{equation}
\frac{1}{t}\mathcal{Y}_{n,\beta }^{(v)}(x;k,a,b)=\frac{1}{n+1}\mathcal{Y}%
_{n+1,\beta }^{(v)}(x;k,a,b).  \label{Re-1}
\end{equation}
\end{theorem}

\begin{theorem}
The following integral formula holds true:%
\begin{equation*}
\left\langle \frac{e^{tc}-1}{t}\mid \mathcal{Y}_{n,\beta
}^{(v)}(x;k,a,b)\right\rangle =\left\{ 
\begin{array}{cc}
\int_{0}^{c}\mathcal{Y}_{n,\beta }^{(v)}(u;k,a,b)du & (n\in \mathbb{N}) \\ 
0 & (n=0).%
\end{array}%
\right. 
\end{equation*}
\end{theorem}

\begin{proof}
We set%
\begin{equation*}
\left\langle \frac{e^{tc}-1}{t}\mid \mathcal{Y}_{n,\beta
}^{(v)}(x;k,a,b)\right\rangle =\left\langle \left( e^{tc}-1\right) \mid 
\frac{1}{t}\mathcal{Y}_{n,\beta }^{(v)}(x;k,a,b)\right\rangle .
\end{equation*}%
By applying the operator in (\ref{Re-1}) and the functional in (\ref{a5}) to
the above equation, we obtain the desired result.
\end{proof}

We give action of a linear operator $\left( \beta ^{b}e^{t}-a^{b}\right) $
on the polynomial $\mathcal{Y}_{n,\beta }^{(v)}(x;k,a,b)$ by the following
lemma:

\begin{lemma}
\label{Lemma3}The following identity holds true:%
\begin{equation*}
\left( \beta ^{b}e^{t}-a^{b}\right) \mathcal{Y}_{n,\beta
}^{(v)}(x;k,a,b)=2^{1-k}\left( n\right) _{k}\mathcal{Y}_{n-k,\beta
}^{(v-1)}(x;k,a,b).
\end{equation*}
\end{lemma}

\begin{proof}
By (\ref{b1}) and Lemma \ref{Lemma1}, we obtain%
\begin{equation*}
\left( \beta ^{b}e^{t}-a^{b}\right) \mathcal{Y}_{n,\beta
}^{(v)}(x;k,a,b)=\left( \beta ^{b}e^{t}-a^{b}\right) \left( \frac{%
2^{1-k}t^{k}}{\beta ^{b}e^{t}-a^{b}}\right) ^{v}x^{n}.
\end{equation*}%
After some calculations, we get%
\begin{equation*}
\left( \beta ^{b}e^{t}-a^{b}\right) \mathcal{Y}_{n,\beta
}^{(v)}(x;k,a,b)=2^{1-k}t^{k}\mathcal{Y}_{n,\beta }^{(v-1)}(x;k,a,b).
\end{equation*}%
Using Lemma \ref{Lemma2} in this equation, it is easy to obtain the desired
result.
\end{proof}

\begin{theorem}
\label{TeoA2}The following identity holds true:%
\begin{equation*}
\mathcal{Y}_{n,\beta }^{(v)}(x+1;k,a,b)=\frac{\left( n\right) _{k}}{%
2^{k-1}\beta ^{b}}\mathcal{Y}_{n-k,\beta }^{(v-1)}(x;k,a,b)+\left( \frac{a}{%
\beta }\right) ^{b}\mathcal{Y}_{n,\beta }^{(v)}(x;k,a,b).
\end{equation*}
\end{theorem}

\begin{proof}
By applying the following well-known operator 
\begin{equation}
e^{at}p\left( x\right) =p\left( x+a\right) \text{ cf. \cite[p. 14]{Roman}}
\label{a17}
\end{equation}%
to the polynomial $\mathcal{Y}_{n,\beta }^{(v)}(x;k,a,b)$, we obtain%
\begin{equation}
\left( \beta ^{b}e^{t}-a^{b}\right) \mathcal{Y}_{n,\beta
}^{(v)}(x;k,a,b)=\beta ^{b}\mathcal{Y}_{n,\beta }^{(v)}(x+1;k,a,b)-a^{b}%
\mathcal{Y}_{n,\beta }^{(v)}(x;k,a,b).  \label{a16}
\end{equation}%
Combining (\ref{a16}) with Lemma \ref{Lemma3}, we arrive at the desired
result.
\end{proof}

\begin{remark}
Taking $\beta =a=b=k=1$ in Lemma \ref{Lemma3} and Theorem \ref{TeoA2},
respectively, we have%
\begin{equation*}
\left( e^{t}-1\right) B_{n}^{\left( v\right) }\left( x\right)
=nB_{n-1}^{\left( v-1\right) }\left( x\right) \text{ \cite[p. 95, Eq-(4.2.5)]%
{Roman}}
\end{equation*}%
and%
\begin{equation*}
B_{n}^{\left( v\right) }\left( x+1\right) =nB_{n-1}^{\left( v-1\right)
}\left( x\right) +B_{n}^{\left( v\right) }\left( x\right) \text{\ \cite[p.
95, Eq-(4.2.6)]{Roman}}.
\end{equation*}
\end{remark}

\begin{remark}
Substituting $\beta =b=1$ , $k=0$ and $a=-1$ into Lemma \ref{Lemma3} and
Theorem \ref{TeoA2}, respectively, we obtain%
\begin{equation*}
\left( e^{t}+1\right) E_{n}^{\left( v\right) }\left( x\right)
=2E_{n}^{\left( v-1\right) }\left( x\right) \text{ \cite[p. 103]{Roman}}
\end{equation*}%
and%
\begin{equation*}
E_{n}^{\left( v\right) }\left( x+1\right) =2E_{n}^{\left( v-1\right) }\left(
x\right) -E_{n}^{\left( v\right) }\left( x\right) \text{ \cite[p. 103,
Eq-(4.2.11)]{Roman}}.
\end{equation*}
\end{remark}

\begin{remark}
Putting $\beta =b=k=1$ and $a=-1$ in Lemma \ref{Lemma3} and Theorem \ref%
{TeoA2}, respectively, we obtain%
\begin{equation*}
\left( e^{t}+1\right) G_{n}^{\left( v\right) }\left( x\right)
=2nG_{n-1}^{\left( v-1\right) }\left( x\right) \text{ \cite[p. 5, Theorem 7]%
{DereSimsek}}
\end{equation*}%
and%
\begin{equation*}
G_{n}^{\left( v\right) }\left( x+1\right) =2nG_{n-1}^{\left( v-1\right)
}\left( x\right) -G_{n}^{\left( v\right) }\left( x\right) .
\end{equation*}
\end{remark}

\begin{remark}
If we put $a=b=k=1$ in Lemma \ref{Lemma3}, and Theorem \ref{TeoA2},
respectively, we obtain%
\begin{equation*}
\left( \beta e^{t}-1\right) \mathcal{B}_{n}^{\left( v\right) }\left( x,\beta
\right) =n\mathcal{B}_{n-1}^{\left( v-1\right) }\left( x,\beta \right)
\end{equation*}%
and%
\begin{equation*}
\mathcal{B}_{n}^{\left( v\right) }\left( x+1,\beta \right) =\frac{n}{\beta }%
\mathcal{B}_{n-1}^{\left( v-1\right) }\left( x,\beta \right) +\frac{1}{\beta 
}\mathcal{B}_{n}^{\left( v\right) }\left( x,\beta \right) .
\end{equation*}
\end{remark}

\begin{remark}
If we set $b=1$, $k=0,$ and $a=-1$ in Lemma \ref{Lemma3} and Theorem \ref%
{TeoA2}, respectively, we obtain%
\begin{equation*}
\left( \beta e^{t}+1\right) \mathcal{E}_{n}^{\left( v\right) }\left( x,\beta
\right) =2\mathcal{E}_{n}^{\left( v-1\right) }\left( x,\beta \right)
\end{equation*}%
and%
\begin{equation*}
\mathcal{E}_{n}^{\left( v\right) }\left( x+1,\beta \right) =\frac{2}{\beta }%
\mathcal{E}_{n}^{\left( v-1\right) }\left( x,\beta \right) -\frac{1}{\beta }%
\mathcal{E}_{n}^{\left( v\right) }\left( x,\beta \right) .
\end{equation*}
\end{remark}

\begin{remark}
By substituting $b=k=1$, and $a=-1$ into Lemma \ref{Lemma3} and Theorem \ref%
{TeoA2}, respectively, we obtain%
\begin{equation*}
\left( \beta e^{t}+1\right) \mathcal{G}_{n}^{\left( v\right) }\left( x,\beta
\right) =2n\mathcal{G}_{n-1}^{\left( v-1\right) }\left( x,\beta \right)
\end{equation*}%
and%
\begin{equation*}
\mathcal{G}_{n}^{\left( v\right) }\left( x+1,\beta \right) =\frac{2n}{\beta }%
\mathcal{G}_{n-1}^{\left( v-1\right) }\left( x,\beta \right) -\frac{1}{\beta 
}\mathcal{G}_{n}^{\left( v\right) }\left( x,\beta \right) .
\end{equation*}
\end{remark}

We give action of a linear operator $\frac{1}{\left( \beta
^{b}e^{t}-a^{b}\right) }$ on the polynomial $\mathcal{Y}_{n,\beta
}^{(v)}(x;k,a,b)$ by the following lemma:

\begin{lemma}
\label{Lemma4}The following identity holds true:%
\begin{equation*}
\frac{1}{\left( \beta ^{b}e^{t}-a^{b}\right) }\mathcal{Y}_{n,\beta
}^{(v)}(x;k,a,b)=\frac{\mathcal{Y}_{n+k,\beta }^{(v+1)}(x;k,a,b)}{%
2^{1-k}\dprod\limits_{j=1}^{k}\left( n+j\right) }.
\end{equation*}
\end{lemma}

\begin{proof}
By (\ref{b1}) and Lemma \ref{Lemma1}, we obtain%
\begin{equation*}
\frac{1}{\left( \beta ^{b}e^{t}-a^{b}\right) }\mathcal{Y}_{n,\beta
}^{(v)}(x;k,a,b)=\frac{1}{\left( \beta ^{b}e^{t}-a^{b}\right) }\left( \frac{%
2^{1-k}t^{k}}{\beta ^{b}e^{t}-a^{b}}\right) ^{v}x^{n}.
\end{equation*}%
After some calculations in this equation, we get%
\begin{equation*}
\frac{1}{\left( \beta ^{b}e^{t}-a^{b}\right) }\mathcal{Y}_{n,\beta
}^{(v)}(x;k,a,b)=\frac{1}{2^{1-k}}\frac{1}{t^{k}}\mathcal{Y}_{n,\beta
}^{(v+1)}(x;k,a,b).
\end{equation*}%
As a consequence, we obtain the desired result by using (\ref{Re-1}).
\end{proof}

\begin{remark}
Substituting $\beta =a=b=k=1$ into Lemma \ref{Lemma4}, we obtain%
\begin{equation*}
\frac{1}{\left( e^{t}-1\right) }B_{n}^{\left( v\right) }\left( x\right) =%
\frac{1}{n+1}B_{n+1}^{\left( v+1\right) }\left( x\right) .
\end{equation*}
\end{remark}

\begin{remark}
If we put $\beta =b=1$, $k=0$, and $a=-1$ in Lemma \ref{Lemma4}, we obtain%
\begin{equation*}
\frac{1}{\left( e^{t}+1\right) }E_{n}^{\left( v\right) }\left( x\right) =%
\frac{1}{2}E_{n}^{\left( v+1\right) }\left( x\right) .
\end{equation*}
\end{remark}

\begin{remark}
Taking $\beta =b=k=1$ and $a=-1$ in Lemma \ref{Lemma4}, we have%
\begin{equation*}
\frac{1}{\left( e^{t}+1\right) }G_{n}^{\left( v\right) }\left( x\right) =%
\frac{1}{2\left( n+1\right) }G_{n+1}^{\left( v+1\right) }\left( x\right) ,%
\text{\cite[p. 3, Lemma 3]{DereSimsek}.}
\end{equation*}
\end{remark}

\begin{remark}
Putting $a=b=k=1$ in Lemma \ref{Lemma4}, we obtain%
\begin{equation*}
\frac{1}{\left( \beta e^{t}-1\right) }\mathcal{B}_{n}^{\left( v\right)
}\left( x\right) =\frac{1}{n+1}\mathcal{B}_{n+1}^{\left( v+1\right) }\left(
x\right) .
\end{equation*}
\end{remark}

\begin{remark}
Upon substituting $b=1$, $k=0$, and $a=-1$ into Lemma \ref{Lemma4}, we obtain%
\begin{equation*}
\frac{1}{\left( \beta e^{t}+1\right) }\mathcal{E}_{n}^{\left( v\right)
}\left( x\right) =\frac{1}{2}\mathcal{E}_{n}^{\left( v+1\right) }\left(
x\right) .
\end{equation*}
\end{remark}

\begin{remark}
If we set $b=k=1$ and $a=-1$ in Lemma \ref{Lemma4}, we obtain%
\begin{equation*}
\frac{1}{\left( \beta e^{t}+1\right) }\mathcal{G}_{n}^{\left( v\right)
}\left( x\right) =\frac{1}{2\left( n+1\right) }\mathcal{G}_{n+1}^{\left(
v+1\right) }\left( x\right) .
\end{equation*}
\end{remark}

If we use Lemma \ref{Lemma3} and Lemma \ref{Lemma4}, then we arrive at the
following corollary:

\begin{corollary}
The following identity holds true:%
\begin{equation}
\frac{\beta ^{b}e^{t}}{\left( \beta ^{b}e^{t}-a^{b}\right) }\mathcal{Y}%
_{n,\beta }^{(v)}(x;k,a,b)=\mathcal{Y}_{n,\beta }^{(v)}(x;k,a,b)+\frac{%
\mathcal{Y}_{n+k,\beta }^{(v+1)}(x;k,a,b)}{2^{1-k}a^{-b}\dprod%
\limits_{j=1}^{k}\left( n+j\right) }.  \label{RRe2}
\end{equation}
\end{corollary}

\begin{theorem}
\label{Teo-2}The following recurrence relation holds true:%
\begin{eqnarray*}
\mathcal{Y}_{n+1,\beta }^{(v)}\left( x;k,a,b\right) &=&\left( x-v\right) 
\mathcal{Y}_{n,\beta }^{(v)}(x;k,a,b)+\frac{vk\mathcal{Y}_{n+1,\beta
}^{(v)}(x;k,a,b)}{n+1} \\
&&-\frac{va^{b}\mathcal{Y}_{n+k,\beta }^{(v+1)}(x;k,a,b)}{%
2^{1-k}\dprod\limits_{j=1}^{k}\left( n+j\right) }.
\end{eqnarray*}
\end{theorem}

\begin{proof}
Setting%
\begin{equation*}
g\left( t\right) =\left( \frac{\beta ^{b}e^{t}-a^{b}}{2^{1-k}t^{k}}\right)
^{v}
\end{equation*}%
in (\ref{a12}), we obtain%
\begin{equation*}
\mathcal{Y}_{n+1,\beta }^{(v)}(x;k,a,b)=\left( x-v\frac{\beta
^{b}e^{t}-kt^{-1}\beta ^{b}e^{t}+kt^{-1}a^{b}}{\beta ^{b}e^{t}-a^{b}}\right) 
\mathcal{Y}_{n,\beta }^{(v)}(x;k,a,b).
\end{equation*}%
After elementary manipulations in this equation, we get%
\begin{eqnarray*}
\mathcal{Y}_{n+1,\beta }^{(v)}(x;k,a,b) &=&x\mathcal{Y}_{n,\beta
}^{(v)}(x;k,a,b)-\frac{v\beta ^{b}e^{t}}{\left( \beta ^{b}e^{t}-a^{b}\right) 
}\mathcal{Y}_{n,\beta }^{(v)}(x;k,a,b) \\
&&+vk\frac{1}{t}\mathcal{Y}_{n,\beta }^{(v)}(x;k,a,b).
\end{eqnarray*}%
Combining the above result with (\ref{Re-1}) and (\ref{RRe2}) gives the
recurrence relation.
\end{proof}

\begin{remark}
If we set $\beta =a=b=k=1$ in Theorem \ref{Teo-2}, we obtain%
\begin{equation*}
B_{n+1}^{\left( v\right) }\left( x\right) =\left( x-v\right) B_{n}^{\left(
v\right) }\left( x\right) -\frac{v}{n+1}B_{n+1}^{\left( v+1\right) }\left(
x\right) +\frac{v}{n+1}B_{n+1}^{\left( v\right) }\left( x\right) .
\end{equation*}%
N\"{o}rlund \cite[p. 145, Eq-(81)]{Norlund} also proved the following
recurrence relations:%
\begin{equation*}
B_{n}^{\left( v+1\right) }\left( x\right) =\left( 1-\frac{n}{v}\right)
B_{n}^{\left( v\right) }\left( x\right) +\frac{n}{v}(x-v)B_{n-1}^{\left(
v\right) }\left( x\right) .
\end{equation*}
\end{remark}

\begin{remark}
If we put $\beta =b=1$, $k=0$, and $a=-1$ in Theorem \ref{Teo-2}, we have%
\begin{equation*}
E_{n+1}^{\left( v\right) }\left( x\right) =\left( x-v\right) E_{n}^{\left(
v\right) }\left( x\right) +\frac{v}{2}E_{n}^{\left( v+1\right) }\left(
x\right) ,
\end{equation*}%
which is proved by N\"{o}rlund \cite[p. 145, Eq-(82)]{Norlund}.
\end{remark}

\begin{remark}
By substituting $\beta =b=k=1$ and $a=-1$ into Theorem \ref{Teo-2}, we have%
\begin{equation*}
G_{n+1}^{\left( v\right) }\left( x\right) =\left( x-v\right) G_{n}^{\left(
v\right) }\left( x\right) +\frac{v}{n+1}G_{n+1}^{\left( v\right) }\left(
x\right) +\frac{v}{2\left( n+1\right) }G_{n+1}^{\left( v+1\right) }\left(
x\right) ,
\end{equation*}%
which is proved by the authors \cite{DereSimsek}.
\end{remark}

\begin{remark}
If we set $a=b=k=1$ in Theorem \ref{Teo-2}, we obtain recurrence formula of
the Apostol-Bernoulli polynomials:%
\begin{equation*}
\mathcal{B}_{n+1}^{\left( v\right) }\left( x,\beta \right) =\left(
x-v\right) \mathcal{B}_{n}^{\left( v\right) }\left( x,\beta \right) -\frac{v%
}{n+1}\mathcal{B}_{n+1}^{\left( v+1\right) }\left( x,\beta \right) +\frac{v}{%
n+1}\mathcal{B}_{n+1}^{\left( v\right) }\left( x,\beta \right) .
\end{equation*}
\end{remark}

\begin{remark}
Substituting $b=1$, $k=0$, and $a=-1$ into Theorem \ref{Teo-2}, we obtain
recurrence formula of the Apostol-Euler polynomials:%
\begin{equation*}
\mathcal{E}_{n+1}^{\left( v\right) }\left( x,\beta \right) =\left(
x-v\right) \mathcal{E}_{n}^{\left( v\right) }\left( x,\beta \right) +\frac{v%
}{2}\mathcal{E}_{n}^{\left( v+1\right) }\left( x,\beta \right) .
\end{equation*}
\end{remark}

\begin{remark}
Putting $b=k=1$ and $a=-1$ in Theorem \ref{Teo-2}, we obtain recurrence
formula of the Apostol-Genocchi polynomials:%
\begin{equation*}
\mathcal{G}_{n+1}^{\left( v\right) }\left( x,\beta \right) =\left(
x-v\right) \mathcal{G}_{n}^{\left( v\right) }\left( x,\beta \right) +\frac{v%
}{n+1}\mathcal{G}_{n+1}^{\left( v\right) }\left( x,\beta \right) +\frac{v}{%
2\left( n+1\right) }\mathcal{G}_{n+1}^{\left( v+1\right) }\left( x,\beta
\right) .
\end{equation*}
\end{remark}

Carlitz \cite{Carliz} established a generalization of the Raabe-type
multiplication formulas for the Bernoulli and Euler polynomials.
Subsequently, Ozden \cite{Ozden} and Ozden \textit{et al}. \cite%
{ozdensimseksrivstava} presented a unification and generalization of the not
only Carlitz \cite{Carliz} result but also Karande \textit{et al}. \cite%
{Karande}\ result. Their proof of the unification and generalization of the
Raabe-type multiplication formulas is based upon the generating function (%
\ref{b1}) with $v=1$.

Here our demonstration of the multiplication formula in Theorem \ref%
{TheoremMULTIPLE} and \ref{RD-1} based upon umbral calculus method.

The multiplication formulas for Bernoulli polynomials also drive from the
Hurwitz zeta function. Raabe-type multiplication formulas are use in many
branches of Mathematics and Mathematical Physics; for example to construct
the Dedekind sums and the Hardy-Berndt sums. One may need to use these
formulas to construct unification and generalization of this sums.

\begin{theorem}
\label{TheoremMULTIPLE}Let $m$ be a positive integer. The following
multiplication formula (Raabe-type multiplication formula) for the
polynomials $\mathcal{Y}_{n,\beta }^{\left( v\right) }(mx;k,a,b)$ holds true:%
\begin{eqnarray*}
&&\mathcal{Y}_{n,\beta }^{\left( v\right) }(mx;k,a,b) \\
&=&m^{n-kv}a^{bv\left( m-1\right) }\dsum\limits_{u_{1},u_{2},\ldots
,u_{m-1}\geq 0}\binom{v}{u_{1},u_{2},\ldots ,u_{m-1}}\left( \frac{\beta }{a}%
\right) ^{bj}\mathcal{Y}_{n,\beta ^{m}}^{\left( v\right) }(x+\frac{j}{m}%
,k,a^{m},b),
\end{eqnarray*}%
where%
\begin{equation*}
j=u_{1}+2u_{2}+\ldots +\left( m-1\right) u_{m-1}.
\end{equation*}
\end{theorem}

\begin{proof}
If we substitute 
\begin{equation*}
g\left( t\right) =\left( \frac{\beta ^{b}e^{t}-a^{b}}{2^{1-k}t^{k}}\right)
^{v}
\end{equation*}%
into (\ref{a14}), we obtain%
\begin{equation*}
\mathcal{Y}_{n,\beta }^{\left( v\right) }(mx;k,a,b)=m^{n-kv}\frac{1}{%
a^{bv}\left( \frac{\beta ^{b}}{a^{b}}e^{\frac{t}{m}}-1\right) ^{v}}\left(
2^{1-k}t^{k}\right) ^{v}x^{n},
\end{equation*}%
From this equation, we obtain%
\begin{equation*}
\mathcal{Y}_{n,\beta }^{\left( v\right) }(mx;k,a,b)=\frac{%
2^{(1-k)v}m^{n-kv}a^{-bv}t^{kv}}{\left( \left( \frac{\beta }{a}\right)
^{bm}e^{t}-1\right) ^{v}}\dsum\limits_{j=0}^{m-1}\left( \frac{\beta }{a}%
\right) ^{bj}e^{\frac{tj}{m}}x^{n}.
\end{equation*}

Thus by applying (\ref{b1}) and Lemma \ref{Lemma1} in the above equation, we
deduce%
\begin{eqnarray*}
&&\mathcal{Y}_{n,\beta }^{\left( v\right) }(mx;k,a,b) \\
&=&m^{n-kv}a^{bvm-bv}\dsum\limits_{u_{1},u_{2},\ldots ,u_{m-1}\geq 0}\binom{v%
}{u_{1},u_{2},\ldots ,u_{m-1}}\left( \frac{\beta }{a}\right) ^{bj}e^{\frac{tj%
}{m}}\left( \frac{2^{1-k}t^{k}}{\beta ^{bm}e^{t}-a^{bm}}\right) ^{v}x^{n}.
\end{eqnarray*}%
By using (\ref{a17}) with Lemma \ref{Lemma1} in the above equation, we
arrive at the desired result.
\end{proof}

\begin{remark}
Theorem \ref{TheoremMULTIPLE} is also related to some special multiplication
formulas. For example, if we set $a=b=k=1$ into assertion of\ Theorem \ref%
{TheoremMULTIPLE}, we arrive at a multiplication formula of the
Apostol-Bernoulli and Apostol-Euler polynomials of higher order, which was
defined by Luo \cite[Theorem 2.1 and Theorem 3.1]{Luo2} and also Bayad and
Simsek \cite[Corollary 5]{bayadsimsek}.
\end{remark}

\begin{remark}
If we set $\beta =a=b=k=1$ in assertion of\ Theorem \ref{TheoremMULTIPLE},
we arrive at the multiplication formula of the Bernoulli polynomials of
higher order. If we take $\beta =b=1$, $k=0$ and $a=-1$ into assertion of\
Theorem \ref{TheoremMULTIPLE}, we arrive at the multiplication formula of
the Euler polynomials of higher order (see, for details, \cite[Corollary 2.1
and Corollary 3.1 ]{Luo2}).
\end{remark}

If we set $v=1$ in Theorem \ref{TheoremMULTIPLE}, we obtain obtain
multiplication formula (Raabe formula) for the 

\begin{corollary}
\label{Teo-3}Let $m$ be a positive integer. The following multiplication
formula holds true:%
\begin{equation}
\mathcal{Y}_{n,\beta }(mx;k,a,b)=m^{n-k}a^{b\left( m-1\right)
}\dsum\limits_{j=0}^{m-1}\left( \frac{\beta }{a}\right) ^{bj}\mathcal{Y}%
_{n,\beta ^{m}}(x+\frac{j}{m},k,a^{m},b).  \label{RD-1}
\end{equation}
\end{corollary}

\begin{remark}
If we substitute $g\left( t\right) =\frac{\beta ^{b}e^{t}-a^{b}}{2^{1-k}t^{k}%
}$ into (\ref{a14}), we obtain%
\begin{equation*}
\mathcal{Y}_{n,\beta }(mx;k,a,b)=\alpha ^{n-k}a^{bm-b}\frac{\frac{\beta ^{b}%
}{a^{b}}e^{t}-1}{\frac{\beta ^{bm}}{a^{bm}}e^{\frac{t}{m}}-1}\mathcal{Y}%
_{n,\beta ^{m}}(x;k,a^{m},b),
\end{equation*}%
By using (\ref{a17}) into the above equation and after some calculations, we
arrive at proof of (\ref{RD-1}).
\end{remark}

\begin{remark}
In special case, when $\beta =a=b=k=v=1$ into assertion of Theorem \ref%
{TheoremMULTIPLE}, we have%
\begin{equation*}
\mathcal{Y}_{n,1}(mx;1,1,1)=B_{n}\left( mx\right)
=m^{n-1}\dsum\limits_{j=0}^{m-1}B_{n}\left( x+\frac{j}{m}\right) ,
\end{equation*}%
(see, for details, \cite{cangul2008}-\cite{Tempesta}; see also the
references cited in each of these earlier works).
\end{remark}

\begin{remark}
In special case, if we take $\beta =b=v=1$, $k=0$ and $a=-1$ into assertion
of Theorem \ref{TheoremMULTIPLE}, we obtain multiplication formulas for
Euler polynomials:
\end{remark}

If $m=2d+1$, $d\in 
\mathbb{N}
$, then%
\begin{equation*}
E_{n}\left( \left( 2d+1\right) x\right) =\left( 2d+1\right)
^{n}\dsum\limits_{j=0}^{2d}\left( -1\right) ^{j}E_{n}\left( x+\frac{j}{2d+1}%
\right) 
\end{equation*}%
(see, for details, \cite{cangul2008}-\cite{Tempesta}; see also the
references cited in each of these earlier works).

If $m=2d$, $d\in 
\mathbb{N}
$, then%
\begin{eqnarray*}
E_{n}\left( 2dx\right)  &=&\left( 2d\right)
^{n}(-1)\dsum\limits_{j=0}^{2d-1}\left( -1\right) ^{j}\mathcal{Y}%
_{n,1}(x;0,-1,1) \\
&=&-\left( 2d\right) ^{n}\dsum\limits_{j=0}^{2d-1}\left( -1\right) ^{j}\frac{%
2}{t}B_{n}\left( x+\frac{j}{2d}\right) 
\end{eqnarray*}

Using (\ref{Re-1}) in the right side of the above equation, we derive the
following\ result:%
\begin{equation*}
E_{n}\left( 2dx\right) =-\frac{2\left( 2d\right) ^{n}}{n+1}%
\dsum\limits_{j=0}^{2d-1}\left( -1\right) ^{j}B_{n+1}\left( x+\frac{j}{2d}%
\right) 
\end{equation*}%
(see, for details, \cite{cangul2008}-\cite{Tempesta}; see also the
references cited in each of these earlier works).

\begin{remark}
Multiplication formulas of the Genocchi polynomials and the Apostol type
polynomials are the special case of Corollary \ref{Teo-3}.
\end{remark}

The next theorem gives us the polynomials $\mathcal{Y}_{n,\beta }^{\left(
v\right) }(x;k,a,b)$ are in terms of the Apostol type polynomials.

\begin{theorem}
Each of the following relationships holds true: 
\begin{equation}
\mathcal{Y}_{n,\beta }^{\left( v\right) }(x;k,a,b)=\frac{\dprod%
\limits_{y=0}^{v\left( k-1\right) -1}\left( n-y\right) }{2^{v\left(
k-1\right) }a^{vb}}\mathcal{B}_{n-v\left( k-1\right) }^{\left( v\right)
}\left( x,\left( \frac{\beta }{a}\right) ^{b}\right) ,  \label{RD-2}
\end{equation}%
\begin{equation}
\mathcal{Y}_{n,\beta }^{\left( v\right) }(x;k,a,b)=-\frac{%
\dprod\limits_{y=0}^{kv-1}\left( n-y\right) }{2^{kv}a^{bv}}\mathcal{E}%
_{n-kv}^{\left( v\right) }\left( x,-\left( \frac{\beta }{a}\right)
^{b}\right) ,  \label{RD-3}
\end{equation}%
and%
\begin{equation}
\mathcal{Y}_{n,\beta }^{\left( v\right) }(x;k,a,b)=-\frac{%
\dprod\limits_{y=0}^{v\left( k-1\right) -1}\left( n-y\right) }{2^{kv}a^{bv}}%
\mathcal{G}_{n-v\left( k-1\right) }^{\left( v\right) }\left( x,-\left( \frac{%
\beta }{a}\right) ^{b}\right) .  \label{RD-4}
\end{equation}
\end{theorem}

\begin{proof}[Proof of (\protect\ref{RD-2})]
\begin{equation*}
\mathcal{Y}_{n,\beta }^{(v)}(x;k,a,b)=\frac{2^{v\left( 1-k\right) }}{a^{vb}}%
t^{v\left( k-1\right) }\mathcal{B}_{n}^{\left( v\right) }\left( x,\left( 
\frac{\beta }{a}\right) ^{b}\right) ,
\end{equation*}%
where%
\begin{equation*}
\mathcal{B}_{n}^{\left( v\right) }\left( x,\left( \frac{\beta }{a}\right)
^{b}\right) =\left( \frac{t}{\left( \frac{\beta }{a}\right) ^{b}e^{t}-1}%
\right) ^{v}x^{n}.
\end{equation*}%
Applying the derivative operator in (\ref{a11}) $v(k-1)$ times in this
equation,$\ $we complete the proof of the (\ref{RD-2}).
\end{proof}

\begin{remark}
We note that proofs of (\ref{RD-3}) and (\ref{RD-4}) are similar to that of (%
\ref{RD-2}). Thus we omit them.
\end{remark}

We now ready to give a relation between the Stirling numbers of the first
kind and the polynomials $\mathcal{Y}_{n,\beta }^{\left( v\right) }(x;k,a,b)$%
.

\begin{theorem}
The following relation holds true:%
\begin{equation*}
\mathcal{Y}_{n,\beta }^{\left( v\right)
}(x;k,a,b)=\dsum\limits_{j=0}^{v}\dsum\limits_{l=0}^{j}\dsum%
\limits_{h=0}^{kl}\left( -1\right) ^{j-l}\binom{v}{j}\binom{j}{l}\frac{%
s\left( kl,h\right) \mathcal{Y}_{n+k\left( v-l\right) ,\beta }^{\left(
v\right) }(x,k,a,b)}{2^{\left( k-1\right) \left( l+v\right)
}n^{-h}\dprod\limits_{j=1}^{kv}\left( n-kl+j\right) },
\end{equation*}%
where $s\left( kl,h\right) $ is the Stirling numbers of the first kind.
\end{theorem}

\begin{proof}
By using Lemma \ref{Lemma1}, we find%
\begin{equation*}
\mathcal{Y}_{n,\beta }^{\left( v\right) }(x;k,a,b)=\dsum\limits_{j=0}^{v}%
\binom{v}{j}\frac{1}{\left( \beta ^{b}e^{t}-a^{b}\right) ^{v}}%
\dsum\limits_{l=0}^{j}\binom{j}{l}\left( -1\right) ^{j-l}\left(
2^{1-k}t^{k}\right) ^{l}x^{n}.
\end{equation*}

By using derivative operator (\ref{a11}) into the above equation, we obtain%
\begin{equation}
\mathcal{Y}_{n,\beta }^{\left( v\right)
}(x;k,a,b)=\dsum\limits_{j=0}^{v}\dsum\limits_{l=0}^{j}\binom{v}{j}\binom{j}{%
l}\frac{1}{\left( \beta ^{b}e^{t}-a^{b}\right) ^{v}}2^{l-kl}\left( -1\right)
^{j-l}\left( n\right) _{kl}x^{n-kl}.  \label{b3}
\end{equation}

Also by using\ Lemma \ref{Lemma1} with 
\begin{equation*}
\left( y\right) _{k}=\dsum\limits_{m=0}^{k}s\left( k,m\right) y^{m}\text{ 
\cite[p. 57]{Roman}}
\end{equation*}%
where $s\left( k,m\right) $ denotes the Stirling numbers of first kind, in
the above equation, we arrive at the desired result.
\end{proof}

We set%
\begin{equation}
\frac{1}{\left( \beta ^{b}e^{t}-a^{b}\right) ^{v}}=-\frac{1}{a^{vb}}%
\dsum\limits_{y=0}^{\infty }\binom{v+y-1}{y}\left( e^{\frac{t}{b}}\frac{%
\beta }{a}\right) ^{yb}.  \label{rr7}
\end{equation}

Observe that further specalizing to $v=1$ the right hand side of equation (%
\ref{rr7}) yields the geometric series expansion.

By substituting (\ref{rr7}) into (\ref{b3}), and after some algebraic
manipulations, we arrive at the following corollary:

\begin{corollary}
The following relation holds true:%
\begin{equation*}
\mathcal{Y}_{n,\beta }^{\left( v\right)
}(x;k,a,b)=\dsum\limits_{y=0}^{\infty
}\dsum\limits_{j=0}^{v}\dsum\limits_{l=0}^{j}\dsum\limits_{h=0}^{kl}\left(
-1\right) ^{j-l}\binom{v+y-1}{y}\binom{v}{j,l,j-l}\frac{s\left( kl,h\right)
\left( x+y\right) ^{h-kl}}{2^{l\left( k-1\right) }a^{b\left( v+y\right)
}\beta ^{-yb}n^{-h}},
\end{equation*}%
where $s\left( kl,h\right) $ is the Stirling numbers of the first kind.
\end{corollary}

\section*{Acknowledgements}

The present investigation was supported, in part, by the \textit{Scientific
Research Project Administration of Akdeniz University}.


\bigskip 

\begin{thebibliography}{99}
\bibitem{bayadsimsek} A. Bayad and Y. Simsek, \textit{Values of twisted
Barnes zeta functions at negative integers}, preprint.

\bibitem{bayad2} A. Bayad and T. Kim, \textit{Identities for the Bernoulli,
the Euler and the Genocchi numbers and polynomials}, Adv. Stud. Contemp.
Math. 20 (2010), pp. 247-253.

\bibitem{blasiack} P. Blasiak, G. Dattoli, A. Horzela and K. A. Penson, 
\textit{Representations of monomiality principle with Sheffer-type
polynomials and boson normal ordering}, Phys. Lett. A 352 (2006), pp. 7-12.

\bibitem{cangul2008} I. N. Cangul, H. Ozden and Y. Simsek, \textit{%
Generating functions of the }$(h,q)$\textit{\ extension of twisted Euler
polynomials and numbers}, Acta Math. Hungar. 120 (2008), pp. 281-299.

\bibitem{cangul2} I. N. Cangul, H. Ozden and Y. Simsek, \textit{A new
approach to }$q$\textit{-Genocchi numbers and their interpolation functions, 
}Nonlinear Anal-Theor. 71 (2009), pp. e793-e799.

\bibitem{Carliz} L. Carlitz, \textit{Some generalized multiplication
formulae for the Bernoulli polynomials and related functions}, Mh. Math. 66
(1962), pp. 1-8.

\bibitem{choi} J. Choi, P. J. Anderson and H. M. Srivastava, \textit{Some }$q
$\textit{-extensions of the Apostol-Bernoulli and the Apostol-Euler
polynomials of order }$n$\textit{, and the multiple Hurwitz zeta function},
Appl. Math. Comput. 199 (2008), pp. 723-737.

\bibitem{choi2009a} J. Choi, P. J. Anderson and H. M. Srivastava, \textit{%
Carlitz's }$q$\textit{-Bernoulli and }$q$\textit{-Euler polynomials and a
class of }$q$\textit{-Hurwitz zeta functions}, Appl. Math. Comput. 215
(2009), pp. 1185-1208.

\bibitem{choi2009b} J. Choi and H. M. Srivastava, \textit{Some applications
of the Gamma and polygamma functions involving convolutions of the Rayleigh
functions, multiple Euler sums and log-sine integrals}, Math. Nachr. 282
(2009), pp. 1709-1723.

\bibitem{Dattoli} G. Dattoli, M. Migliorati and H. M. Srivastava, \textit{%
Sheffer polynomials, monomiality principle, algebraic methods and the theory
of classical polynomials}, Math. Comput. Modelling 45 (2007), pp. 1033-1041.

\bibitem{DereSimsek} R. Dere and Y. Simsek, \textit{Genocchi polynomials
associated with the Umbral algebra}, In press, accepted manuscript, Appl.
Math. Comput. 217 (2011).

\bibitem{KDilcher} K. Dilcher, \textit{Asymptotic behaviour of Bernoulli,
Euler, and generalized Bernoulli polynomials}, J. Approx. Theory 49 (1987),
pp. 321-330.

\bibitem{GJS} M. Garg, K. Jain and H. M. Srivastava, \textit{Some
relationships between the generalized Apostol-Bernoulli polynomials and
Hurwitz-Lerch Zeta functions}, Integral Transform. Spec. Funct. 17 (2006),
pp. 803-815.

\bibitem{jang1} L.-C. Jang, B. Lee and T. Kim, On the twisted $q$-analogs of
the generalized Euler numbers and polynomials of higher order, Adv.
Difference Equ. (Article ID 875098), 2010 Volumes 11 pages (2010).

\bibitem{Karande} B.K. Karande and N.K. Thakare, \textit{On the unification
of Bernoulli and Bernoulli polynomials}, Indian J. Pure Appl. Math. 6(1)
(1975), pp. 98-107.

\bibitem{Kimadcm} T. Kim, \textit{Some identities for the Bernoulli, the
Euler and the Genocchi numbers and polynomials}, Adv. Stud. Contemp. Math.
20 (2010), pp. 23-28.

\bibitem{simsekKimDkim} T. Kim, S.-H. Rim, Y. Simsek, and D. Kim, \textit{On
the analogs of Bernoulli and Euler numbers, related identities and zeta and }%
$L$\textit{-functions}, J. Korean Math. Soc. 45 (2008), pp. 435-453.

\bibitem{Jang} Y.-H. Kim, W. Kim and L.-C. Jang, \textit{On the }$q$\textit{%
-Extension of Apostol-Euler Numbers and polynomials, }Abstract Appl. Anal.
(Article ID 296159), 2008 Volumes 10 pages (2008).

\bibitem{Luo2} Q.-M. Luo, \textit{The multiplication formulas for the
Apostol-Bernoulli and Apostol-Euler polynomials of higher order}, Integral
Transforms Spec. Funct. 20 (2009), pp. 377-391.

\bibitem{Luo3} Q.-M. Luo, $q$\textit{-Extensions for the Apostol-Genocchi
Polynomial}, Gen. Math. 17 (2009), pp. 113-125.

\bibitem{LuoSrivastava} Q.-M. Luo and H. M. Srivastava, \textit{Some
generalizations of the Apostol-Bernoulli and Apostol-Euler polynomials}, J.
Math. Anal. Appl. 308(1) (2005), pp. 290-302.

\bibitem{sir3} Q.-M. Luo and H. M. Srivastava, $q$\textit{-Extensions of
some relationships between the Bernoulli And Euler polynomials}, Taiwanese
J. Math. 15 (2011), pp. 241-257.

\bibitem{sir4} Q.-M. Luo and H. M. Srivastava, \textit{Some generalizations
of the Apostol-Genocchi polynomials and the Stirling numbers of the second
kind}, Appl. Math. Comput. 217(2011), pp. 5702-5728.

\bibitem{Maldonado} M. Maldonado, J. Prada and M. J. Senosiain, \textit{%
Basic Appell Sequences}, Taiwan J. Math. 11 (2007), pp. 1045-1055.

\bibitem{Norlund} N. E. N\"{o}rlund, \textit{Vorlesungen uber
Differenzenrechnung}, Springer 1924.

\bibitem{Ozden} H. Ozden, \textit{Unification of generating function of the
Bernoulli, Euler and Genocchi numbers and polynomials, Numerical Analysis
and Applied Mathematics}, AIP Conf. Proc. 1281, (2010), pp. 1125-1227.

\bibitem{OzdenICNAAM2011} H. Ozden, \textit{Generating functions of the
unified representation of the Bernoulli, Euler and Genocchi polynomials of
higher order} to appear in Numerical Analysis and Applied Mathematics, AIP
Conf. Proc. (2011).

\bibitem{ozdenSimsekAML} H. Ozden and Y. Simsek, \textit{A new extension\ of 
}$q$\textit{-Euler numbers and polynomials related to their interpolation
functions},\ Appl. Math. Lett. 21 (2008), pp. 934-939.

\bibitem{ozdenCangulSimsek} H. Ozden, I.N. Cangul and Y. Simsek, \textit{%
Multivariate interpolation functions of higher-order }$q$\textit{-Euler
numbers and their applications}, Abstract Appl. Anal. (Article ID 390857),
2008 Volumes 16 pages (2008).

\bibitem{ozdensimseksrivstava} H. Ozden, Y. Simsek and H. M. Srivastava, 
\textit{A unified presentation of the generating functions of the
generalized Bernoulli, Euler and Genocchi polynomials}, Comput. Math. Appl.
60 (2010), pp. 2779-2787.

\bibitem{Rim} S.-H. Rim, Y.-H. Kim, B. Lee, and T. Kim, \textit{Some
identities of the generalized twisted Bernoulli numbers and polynomials of
higher order}, J. Comput. Anal. Appl. 12 (2010), pp. 695-702.

\bibitem{Roman} S. Roman, \textit{The Umbral Calculus}, Dover Publ. Inc. New
York, 2005.

\bibitem{SIMSEKjnt-2005} Y. Simsek, $q$\textit{-Analogue of the twisted }$l$%
\textit{-series and }$q$\textit{-twisted Euler numbers}, J. Number Theory
110 (2005), pp. 267-278.

\bibitem{simtwistedhqbernoullijmaa} Y. Simsek, \textit{Twisted }$(h,q)$%
\textit{-Bernoulli numbers and polynomials related to twisted }$(h,q)$%
\textit{-zeta function and }$L$\textit{-function}, J. Math. Anal. Appl. 324
(2006), pp. 790-804.

\bibitem{simsekASCM2008} Y. Simsek, \textit{Generating functions of the
twisted Bernoulli numbers and polynomials associated with their
interpolation functions}, Adv. Stud. Contemp. Math. 16 (2008), pp. 251-278.

\bibitem{simsek2009} Y. Simsek,\textit{\ }$q$\textit{-Hardy-Berndt type sums
associated with }$q$\textit{-Genocchi type zeta and }$q$\textit{-}$l$\textit{%
-functions}, Nonlinear Anal-Theor. 71 (2009), pp. e377-e395.

\bibitem{SimsekSrivastava} Y. Simsek and H. M. Srivastava, \textit{A family
of }$p$\textit{-adic twisted interpolation functions associated with the
modified Bernoulli numbers}, Appl. Math. Comput. 216 (2010), pp. 2976-2987.

\bibitem{HMS} H. M. Srivastava, \textit{Some formulas for the Bernoulli and
Euler polynomials at rational arguments}, Math. Proc. Cambridge Philos. Soc.
129 (2000), pp. 77-84.

\bibitem{pinter} H. M. Srivastava and A. Pinter, Remarks on some
relationships between the Bernoulli and Euler polynomials, Appl. Math. Lett.
17 (2004) pp. 375-380.

\bibitem{sir1} H. M. Srivastava, M. Garg and S. Choudhary, \textit{A new
generalization of the Bernoulli and related polynomials}, Russian J. Math.
Phys. 17 (2010), pp. 251-261.

\bibitem{sir2} H. M. Srivastava, M. Garg and S. Choudhary, \textit{Some new
families of generalized Euler and Genocchi polynomials}, Taiwanese J. Math.
15 (2011), pp. 283-305.

\bibitem{sirKIMsim} H. M. Srivastava, T. Kim and Y. Simsek, $q$\textit{%
-Bernoulli numbers and polynomials associated with multiple }$q$\textit{%
-zeta functions and basic }$L$\textit{-series}, Russian J. Math. Phys. 12
(2005), pp. 241-268.

\bibitem{Tempesta} P. Tempesta, \textit{On Appell sequences of polynomials
of Bernoulli and Euler type}, J. Math. Anal. Appl. 341 (2008), pp. 1295-1310.
\end{thebibliography}
\end{document}